\documentclass[11pt,oneside,psamsfonts,final]{amsart}

\usepackage[body={5.45in,9.3in}]{geometry}

\usepackage[mathscr]{eucal}
\usepackage{amsfonts,amssymb,amsmath,amsthm,graphicx,
mathptm,
times}

\usepackage{hyperref}
\usepackage[all]{hypcap} 

\usepackage[oldenum]{paralist} 

\usepackage{soul} 
\sodef\TopTiFo{}{0.03em}{.5em}{2em plus 0.1em minus 0.1em}

\usepackage{hyperref}

\makeatletter

\def\@settitle{\begin{center}%
\uppercasenonmath\@title{
  \large\textrm{\@title}}
  \end{center}%
}
\copyrightinfo{}{}

\def\@setauthors{%
  \begingroup
  \trivlist
  \centering\footnotesize \@topsep40\p@\relax
  \advance\@topsep by -\baselineskip
  \item\relax
  \andify\authors
\textsc{\authors}%
  \endtrivlist
{\hfill \scriptsize 
        (\textit{Received $8$ August $1989$}; 
         \textit{in revised form $7$ January $1991$}; \hfill}
\par
{\hfill \scriptsize 
        \textit{annotated and \textup{arXiv}ed by the author %
               $13$ November $2004$}) \hfill}
  \endgroup
}
\def\@seccntformat#1{\bfseries
       \protect\footnotesize{$\mathbf{\S}$\csname the#1\endcsname.\quad}}
\newdimen\@bls                              
\@bls=\baselineskip                         
\def\section{\@startsection{section}{1}{\z@}{1.5\@bls\@plus .4\@bls 
             \@minus .1\@bls}{\@bls}{\centering\footnotesize\bfseries}}

\def\@captionheadfont{\normalfont\footnotesize}
\def\@captionfont{\normalfont\footnotesize}
\newtheoremstyle{TopologyThm}{\baselineskip}{\baselineskip}{\itshape}%
  {\parindent}{\fontshape{sc}\selectfont}{.}{ }{}
\newtheoremstyle{TopologyRmk}{\baselineskip}{\baselineskip}{}%
  {\parindent}{\fontshape{it}\selectfont}{.}{ }{}

\theoremstyle{TopologyThm}

\newtheorem{theorem}{T{\footnotesize{heorem}}}[section]
\newtheorem{lemma}[theorem]{L{\footnotesize{emma}}}
\newtheorem*{conjecture}{C{\footnotesize{onjecture}}}

\theoremstyle{TopologyRmk}
\newtheorem*{remark}{Remark}

\def\footnoterule{\kern-.4\p@
        \hrule\@width 5.5pc\kern11\p@\kern-\footnotesep}
\def\@setthanks{\def\thanks##1{\kern-\parindent$\dagger$##1\@addpunct.}\thankses}

\def\article@logo{%
  \set@logo{%
    \textup{Hyper{\TeX}ed arXival version prepared 
    and updated November}\\
    \textup{2004. Marginal notes refer to the 
    \protect\hyperlink{addenda}{Addenda}.
    Original publication:}\\
    \textup{\textit{Topology}, Vol.\ 31, No.\ 2, pp.\ 231--237, 1992.}
  }%
}

\DeclareRobustCommand{\appnote}[1]{%
\marginpar{\quad\label{appnote:#1}%
\hypertarget{appnote:#1}{\large\textcircled{\small\ref{#1}}}}}

\def\qedsymbol{{\ensuremath\Box}}
\newenvironment{unproved}[1]{\def\QQQQ{\string#1}\begin{\QQQQ}}%
                                  {\hfill\qedsymbol\end{\QQQQ}}

\newcommand\Bd{\partial} 
\newcommand{\deffont}{\textit}        
\newcommand{\bydef}[1]{\deffont{#1}}  
\newcommand{\Suchthat}{\negthinspace\colon} 
\newcommand{\C}{\mathbb{C}}
\def\s#1{{{\scriptstyle\sigma\mathstrut}_{\scriptscriptstyle{#1}}}}
\newcommand{\e}{\varepsilon}
\newcommand{\sub}{\subset}
\newcommand{\R}{\mathbb{R}}
\renewcommand\P{\mathbb{P}}

\newcommand{\shout}[1]{\emph{#1}}     
\newcommand{\emband}[2]{{\scriptstyle\sigma}_{{\scriptscriptstyle{#1,\mkern3mu #2}}}} 
\newcommand\brep[1]{\mathbf{#1}}
\newcommand\h[2]{h^{\scriptscriptstyle(#1)}_{#2}} 
\newcommand\CH[2]{H^{\scriptscriptstyle(#1)}_{#2}} 
\newcommand\isdefinedas{\mathrel{:=}}
\newcommand\skel[2]{#1^{{\scriptscriptstyle(#2)}}} 
\newcommand\Int{\operatorname{Int}} 

\newcommand\overstrike[2]{{\setbox0\hbox{$#2$}\hbox to \wd0{\hss
                         $#1$\hss}\kern-\wd0\box0}}

\newcommand\connsum{\mathop{\hbox{\overstrike%
         {\hbox{\raise-.25ex\hbox{/\kern -.08em/}}}=}}}
\newcommand\Connsum{\mathop{\displaystyle\hbox{\overstrike%
         {\hbox{\raise-.25ex\hbox{/\kern -.08em/}}}=}}}

\begin{document}

\bibliographystyle{amsplain}
\title[Quasipositive Seifert surfaces]%
{\TopTiFo{A characterization of quasipositive Seifert}\\
\TopTiFo{surfaces (Constructions of quasipositive knots and}\\
\TopTiFo{links, III)}}
\author[Lee Rudolph]%
{\large{L}\small{ee} \large{R}\small{udolph}\thinspace$\dagger$}

\thanks{Partially supported by NSF grant 
\href{http://www.nsf.gov/cgi-bin/showaward?award=8801959}
{DMS-8801959}.}

\def\ufootnote#1{\let\savedthfn\thefootnote\let\thefootnote\relax
\footnote{#1}\let\footnotemark\savedthfn\addtocounter{footnote}{-1}}

\def\dots{\mathinner{\mkern2mu\ldotp\mkern2mu\ldotp\mkern2mu\ldotp\mkern2mu}}
\def\dotsm{\mathinner{\mkern2mu\cdotp\mkern2mu\cdotp\mkern2mu\cdotp\mkern2mu}}

\maketitle
\makeatother

\setcounter{section}{-1}

\pagestyle{myheadings}
\markright{\protect\small Lee Rudolph}

\section{INTRODUCTION; STATEMENT OF RESULTS}

\noindent\textsc{\large{L}\small{et}} $T_{n,n}\sub S^3$  be a fiber surface 
for the torus link $O\{n,n\}$; say, to be concrete,\break
$T_{n,n} = \{(z,w){\thinspace\scriptstyle\in\thinspace} 
S^3\Suchthat z^n+w^n\ge 0\} \sub 
S^3 = \{(z,w)\Suchthat |z|^2 +|w|^2=1\}$.

\begin{unproved}{CharThm}
A Seifert surface $S$ is quasipositive if and only 
if\textup{,} for some~$n$\textup{,} $S$ is ambient 
isotopic to a full subsurface of $T_{n,n}$. 
\end{unproved}

This should be contrasted with the following.

\begin{unproved}{LyonThm}
Any Seifert surface is ambient isotopic to a \textup(full\textup) 
subsurface of the fiber surface of $O\{n,n\}\connsum O\{n,n\}$ for some $n$.
\end{unproved}

Here, a \bydef{surface} is smooth, compact, oriented, and has 
no component with empty
boundary.  A \bydef{Seifert surface}
is a surface embedded in $S^3$.  A subsurface 
$S$ of a surface $T$ is \bydef{full} 
if each simple 
closed curve on $S$ that bounds a disk on $T$ already bounds a 
disk on $S$.  The definition of quasipositivity 
is recalled in \S\ref{braids}, after a review of braided surfaces.

The ``only if'' statement of the Characterization Theorem is proved in 
\S\ref{second section}.  Some results  
about graphs on braided surfaces (stated, with one 
eye on other applications \cite{Rudolph1990b}, in  
somewhat more generality 
than needed here) are obtained in \S\ref{graphs}, 
and used in \S\ref{fourth section} to prove the 
``if'' statement of the Characterization Theorem. A conjectural 
extension to ribbon surfaces   
in the $4$-disk is discussed in \S\ref{generalization}.

\begin{remark}
The relation ``$S$ is a full subsurface of $T$'' is transitive, so 
the Characterization Theorem has the interesting corollary (applied in 
\cite{Rudolph1990a} and 
\cite{Rudolph1989}) that \shout{any full subsurface of any quasipositive 
Seifert surface is quasipositive}.  
\end{remark}

\section{BRAIDS, BRAIDED SURFACES, AND QUASIPOSITIVITY}\label{braids}

Let $\s1,\dots,\s{n-1}$ be the usual generators 
of the $n-$string braid group $B_n$.  For\break  
$1\le i<j\le n$, let  $\emband{i}{j}\isdefinedas 
(\s{i}\dotsm \s{j-2})\s{j-1}(\s{i}\dotsm \s{j-2})^%
{\scriptscriptstyle{-1}}$.  
A \bydef{positive} (resp., \bydef{negative})\break  
\bydef{embedded band} in $B_n$ is any 
$\emband{i}{j}$ (resp., $\emband{i}{j}^{\scriptscriptstyle{-1}}$).  
An \bydef{embedded band representation} in $B_n$, of length
$k$, is a $k$-tuple $\brep{b}=(b(1),\dots,b(k))$ of embedded bands.  
If each $b(t)$ is positive, then $\brep{b}$ is\linebreak[4]

\pagebreak[4]\noindent
\bydef{quasipositive}.  If each $b(t)$ is some 
usual generator $\s{i}=\emband{i}{i+1}$
or inverse $\s{i}^{\scriptscriptstyle{-1}}$, then $\brep{b}$ is 
a \bydef{braidword}. (A quasipositive braidword 
is customarily simply called 
\bydef{positive}.)

\begin{figure}
\centering
\includegraphics[width=2.0in]{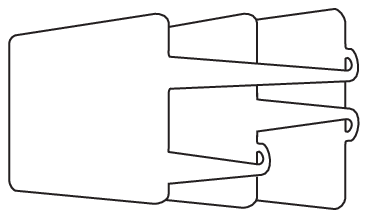}
\caption{A braided Seifert surface $S(\brep{b})$,
$\brep{b}=(\emband{1}{2}^{\protect\mathstrut},
\emband{2}{3}^{\scriptscriptstyle{-1}\protect\mathstrut},
\emband{1}{3}^{\protect\mathstrut})$.
\label{Figure 1}}
\end{figure}

Given $\brep{b}$, a \bydef{braided surface}
$S(\brep{b})$ is constructed as in Fig.~\ref{Figure 1} 
(cf.~\cite{Rudolph1983b}\hyperlink{convention}{$\ddagger$}\ufootnote{%
  \kern-\parindent\hypertarget{convention}{$\ddagger$}\thinspace 
  The orientation conventions in 
  \cite{Rudolph1983b}--\cite{Rudolph1984} are opposed to those of
  \cite{Rudolph1990a}, \cite{Rudolph1990b}, and the present paper.}).  
More precisely, if 
$b(t)=\emband{i(t)}{j(t)}^{\e(t)} {\thinspace\scriptstyle\in\thinspace} B_n$, 
$1\le i(t) < j(t)\le n$, $\e(t)=\pm1$, $1\le t\le k$, then 
$S(\brep{b})\sub\R^3\sub S^3$ is 
a Seifert surface with handle decomposition 
$\bigcup\limits_{s=1}^n \h0s \cup
\bigcup\limits_{t=1}^k \h1t $
such that, for appropriate 
rectangular coordinates 
$(x,y,z)$ on $\R^3$,
\begin{enumerate}
\item\label{zero-handles}
$\h0s$ is contained in the vertical half plane
$\{(s,y,z)\Suchthat y\ge0, z{\thinspace\scriptstyle\in\thinspace}\R\}$,
\item
$\h1t$ is contained in the box 
$\{(x,y,z)\Suchthat i(t)\le x\le j(t), y\le 0, t-1\le z\le t \}$,
\item
$\h1t$ joins $\h0{i(t)}$ and $\h0{j(t)}$ with a single half-twist 
of sign $\e(t)$.
\end{enumerate}

According to \cite{Rudolph1983b}, 
this construction has a converse: if $S$ is any Seifert 
surface, then there exists $\brep{b}$ (highly nonunique!) 
such that $S$ is ambient isotopic to $S(\brep{b})$.
$S$ is \bydef{quasipositive} if some such $\brep{b}$ is quasipositive.

Note that (after enlarging its $0$-handles in their 
half-planes as necessary) the link $\Bd S(\brep{b})$~is 
presented, with respect to a horizontal braid axis (omitted 
from the figure) as the closed %
braid $\widehat\beta(\brep{b})$, 
where $\beta(b)\isdefinedas b(1)\dotsm b(k)$.

\section{}\label{second section}
The ``only if'' statement in the Characterization Theorem 
follows immediately from \ref{lemma 2.1} 
and~\ref{lemma 2.2}, plus the observation that the relation 
``$S$ is a full subsurface of $T$'' is transitive.

\begin{lemma}\label{lemma 2.1}
If $\brep{p}$ is a positive braidword 
in $B_m$\textup{,} then\textup{,}\thinspace{for} some $n$\textup{,} $S(p)$ 
is ambient isotopic to a full subsurface of $T_{n,n}$.
\end{lemma}

\begin{proof}
Let $\mathbf{\nabla}_n$  be the positive braidword of length $n(n-1)$ 
in $B_n$ with $\mathbf{\nabla}(s)=\s{d}$ if 
$s=(n-1)c+d$, $0\le c\le n-1$, $1\le d\le n-1$.  
Then $\widehat\beta(\mathbf{\nabla}_n)$ is a 
torus link $O\{n,n\}$ of~type 
$(n,n)$, and $S(\mathbf{\nabla}_n)$ 
is a fiber surface $T_{n,n}$ (cf. \cite{Stallings}).  
Let $n$ be 
the greater of $m$ and the length~of~the braidword $\brep{p}$; 
make the usual identification of $B_m$  with a subgroup of $B_n$.  
By lavishly 
inserting letters into $\brep{p}$, 
pad it out to $\mathbf{\nabla}_n$.  
Now $S(\brep{p})$ may be constructed as a (manifestly full) 
subsurface of $S(\mathbf{\nabla}_n)$. 
\end{proof}

\begin{lemma}\label{lemma 2.2}
If $S$ is a quasipositive Seifert surface, then $S$ is ambient 
isotopic to a full subsurface of $S(\brep{p})$ for some positive 
braidword $\brep{p}$ in some $B_n$ .
\end{lemma}
\pagebreak[4]
\markright{\protect\small QUASIPOSITIVE SEIFERT SURFACES}
\begin{figure}
\centering
\includegraphics[width=4.75in]{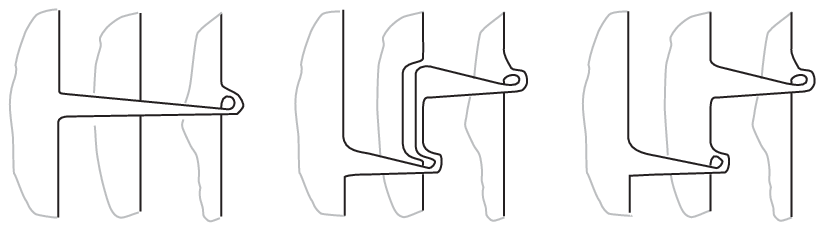}
\caption{\label{Figure 2}}
\end{figure}

\begin{proof}
We may assume that $S=S(\brep{b})$, where 
$\brep{b}=(\emband{i(1)}{j(1)},\dots,\emband{i(k)}{j(k)})$ is 
a quasipositive embedded band representation in some $B_n$.  
By inspection (cf. the $1$-handle pictured~in 
Fig.~\ref{Figure 2}, 
$S(\brep{b})$ is ambient isotopic to a 
(full) subsurface of $S(\brep{p})$ where 
$\brep{p}$ is the positive 
braidword 
$(\s{i(1)},\s{i(1)+1},\dots,\s{j(1)-1},
\s{i(2)},\dots,\s{j(2)-1},\dots,\s{i(k)},\dots,\s{j(k)-1})$. 
\end{proof}

\section{GRAPHS ON BRAIDED SURFACES}\label{graphs}
A \bydef{graph} is a polyhedron $G$ of dimension $\le 1$.  
If $g{\thinspace\scriptstyle\in\thinspace} G$, 
then $g$ is an \bydef{isolated point} (resp., 
\bydef{endpoint}; \bydef{ordinary point}; \bydef{intrinsic vertex}) 
if the link of $g$ in $G$ consists of $0$ (resp., $1$; $2$; $v(g)\ge 3$) 
points.  $G$ is \bydef{trivalent} if $v(g)=3$ for each 
intrinsic vertex $g{\thinspace\scriptstyle\in\thinspace} G$.  

A graph $E$ contained in a vertical half-plane 
$\{(s,y,z)\Suchthat y\ge 0, z{\thinspace\scriptstyle\in\thinspace} \R\}$ 
is a \bydef{comb} if $E$ is the 
union of a vertical interval 
$\{(s,y,z)\Suchthat y=y_0 , z_1\le z\le z_m\}$ and 
$m\ge 2$ horizontal intervals 
$\{(s,y,z)\Suchthat x=s, y_0\ge y\ge 0, z=z_h\}$, 
$y_0>0$, $z_1<\dotsm< z_m$; 
$E$ is trivalent, with $m$ endpoints (on the boundary of the half-plane) 
and $m-2$ intrinsic vertices.  

Let $S$ be a surface in which a graph $G$ is embedded.  
Let $N_S(G)$ be a regular neighborhood of $G$.  
$G$ is \bydef{full} if no simple closed curve $C\sub G$ 
bounds a disk on $S$ (so $G$ is full if and only if $N_S(G)$ is), 
and $G$ is proper if $G\cap\Bd S$ is precisely the set of 
endpoints of $G$.  \mbox{If $N\sub S$ is} 
a subsurface such that 
each component of $N\cap \Bd S$ is an arc, then 
$N$ is a regular neighborhood $N_S(G)$ of some proper trivalent graph $G\sub S$.  

Let $G'$ be another graph embedded in $S$.  
$G$ and $G'$ are \bydef{isotopic} if there exists 
a piecewise-smooth isotopy of $S$ carrying $G$ onto $G'$.  
(It may not be possible to find such an isotopy 
which is smooth near intrinsic vertices.)  $G$ and $G'$ are 
\bydef{Whitehead-equivalent} if $N_S(G)$ and $N_S(G')$ 
are ambient isotopic on S.    

Given a handle decomposition 
$S = \bigcup\limits_{s=1}^n \h0s \cup \bigcup\limits_{t=1}^k \h1t$,
set $\skel{G}{0}\isdefinedas G\cap \bigcup\limits_{s=1}^n \h0s$, 
$\smash[t]{\skel{G}{1}\isdefinedas G\cap \bigcup\limits_{t=1}^k \h1t}$.  
$G$ is \bydef{well-placed} if each point of $\skel{G}{1}$ 
is an ordinary point of $G$, and each     
component of $\skel{G}{1}$ is isotopic to a core arc of some $\h1t$.  
If $S=S(\brep{b})$ is braided and    
$\bigcup\limits_{s=1}^n \h0s \cup \bigcup\limits_{t=1}^k \h1t$
is its given handle decomposition, $G$ is \bydef{combed} 
if it is well-placed and $\skel{G}{0}$     
is a disjoint union of combs and isolated points.

\begin{lemma}\label{lemma 3.1}
Any graph on S is isotopic to a well-placed graph.
\end{lemma}

\begin{proof}
Obvious. 
\end{proof}

\begin{lemma}\label{lemma 3.2}
Any trivalent, full, proper graph on $S(\brep{b})$ is isotopic, by 
an isotopy supported on the $0$-handles, to a combed graph.
\end{lemma}
\pagebreak[4]
\markright{\protect\small Lee Rudolph}
\begin{proof}
Let $G$ satisfy the hypotheses.  By \ref{lemma 3.1}, 
we can assume G is well-placed.  Let 
$S(\brep{b})\cap\{(s,y,z)\Suchthat y\ge 0, 
z{\thinspace\scriptstyle\in\thinspace}\R\} = \h0s$ be a $0$-handle.  
By fullness of $G$ and the Jordan Curve Theorem, 
each component of $G\cap\h0s$  is an isolated point 
or a tree.  By the properness~of~$G$, 
if $e$ is an endpoint of $G\cap\h0s$, 
then $e{\thinspace\scriptstyle\in\thinspace}\Bd\h0s$, 
and $e$ is interior to an attaching arc of some 
$1$-handle (resp., exterior to all the attaching arcs) 
iff $e$ is an ordinary point (resp., an 
endpoint) of $G$.  
Thus after a preliminary isotopy supported 
in a collar of $\Bd\h0s$, we can 
assume that every 
endpoint of $G\cap\h0s$  is an interior point of the 
interval 
$J\isdefinedas S(\brep{b})\cap
\{(s,0,z)\Suchthat z{\thinspace\scriptstyle\in\thinspace}\R\}$.  
For each tree component $E$ of $G\cap\h0s$, 
let $I\sub J$ be the smallest subinterval containing 
all the endpoints of $E$, let $A$ be the subarc of $E$ 
joining the endpoints 
of $I$, and let 
$D(E)$ be the subdisk of $\h0s$ bounded by $I\cup A$.  The set 
of disks $D(E)$~is~partially ordered by inclusion.  
Using trivalence of $G$, we can construct the desired 
isotopy by 
induction over this poset.  
\end{proof}

\begin{lemma}\label{lemma 3.3}
Let $G\sub S(\brep{b})$ be combed.  
If for some $s$\textup{,} $t$\textup{,} 
there is a comb $E\sub G\cap\h0s$ which has two or more 
endpoints in the attaching arc $\h0s\cap\h1t$\textup{,} then 
there is a combed graph $G'$ Whitehead-equivalent to $G$ 
such that  $\skel{G'}{1}$ has fewer components than 
$\skel{G}{1}$.
\end{lemma}

\begin{figure}[h]
\centering
\includegraphics[width=1.875in]{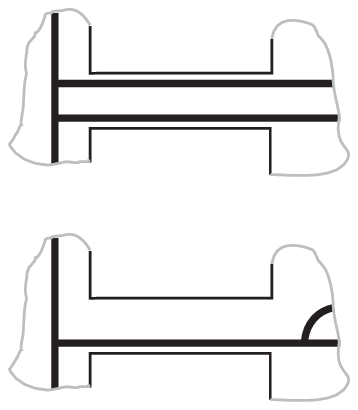}
\caption{\label{Figure 3}}
\end{figure}

\begin{proof}
Let $E=\{(s,y_0,z)\Suchthat z_1\le z\le z_m\}\cup%
\bigcup\limits_{i=1}^m \{(s,y,z)\Suchthat y_0\ge y\ge 0, z=z_i\}$ 
be such  
a comb.  Then, for some $i$, adjacent endpoints 
$(s,0,z_i)$, $(s,0,z_{i+1})$ of $E$ are in 
$\h1t$.  Altering 
$G$ in a neighborhood of 
$\h1t\cup\{(s,y,z)\Suchthat 0\le y\le y_0, z_i\le z\le z_{1+1}\}$, 
we can replace $G$ with 
a Whitehead-equivalent graph $G''$     
which is still trivalent, full, proper, and well-placed, such 
that $\skel{G''}{1}$ has fewer components than $\skel{G}{1}$.  
(Cf.~Fig.~\ref{Figure 3}; the ``Whitehead move'' we use is 
the equivalent, on the level of graphs, of a handle-slide 
on the level of regular neighborhoods.) 
Then we apply \ref{lemma 3.2} to $G''$ to obtain $G'$. 
\end{proof}

\section{}\label{fourth section}
\begin{figure}[t]
\centering
\includegraphics[width=1.625in]{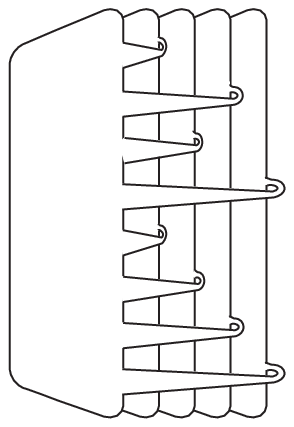}
\caption{\label{Figure 4}}
\end{figure}

At this point, $S(\mathbf\nabla_{n})$ is no longer the best 
choice of a braided surface ambient isotopic to  
$T_{n,n}$.  Instead, for $\nu\isdefinedas (n-1)^2+1$, let 
$\brep{q}$ be the quasipositive band representation~of~length 
$2(\nu-1)$ in $B_\nu$ defined as follows: 
for $1\le s\le \nu-1$, set $q(s)=\emband{1}{\nu-s+1}$; 
for  
$\nu\le s\le 2(\nu-1)$, 
if $s-\nu = (n-1)c+d$, $0\le c, d\le n-2$, 
then set $q(s)=\emband{1}{\nu-c-(n-1)d}$.     
(The case $n=3$ is illustrated in Fig.~\ref{Figure 4}.) 
\pagebreak[4]
\markright{\protect\small QUASIPOSITIVE SEIFERT SURFACES}
\begin{lemma}\label{lemma 4.1}
$S(\brep{q})$ is ambient isotopic to $T_{n,n}$.
\end{lemma}

\begin{proof}
First show that $\widehat\beta(\brep{q})$ is ambient isotopic 
to $O\{n,n\}$ (this is a straightforward 
exercise in braid relations and Markov moves, given that 
$O\{n,n\}=\widehat\beta(\mathbf{\nabla}_n))$.  
Then calculate 
\[
\chi(S(\brep{q}))=\nu-2(\nu-1)=2-\nu=1-(n-1)^2= n-n(n-1)=\chi(T_{n,n}).  
\]
Finally, use the familiar fact that up to ambient isotopy 
a fiber surface is the unique Seifert surface of maximal 
Euler characteristic for its boundary.  
\end{proof}

In addition to its given handle decomposition
$S(\brep{q})=\bigcup\limits_{s=1}^\nu \h0s \cup
\bigcup\limits_{t=1}^{2(\nu-1)} \h1t$, which will 
now be called \bydef{fine}, we need a \bydef{coarse} 
handle decomposition with a single $0$-handle 
$\CH01\isdefinedas\h01$ and $\nu-1$ $1$-handles 
$\CH1s\isdefinedas\h1s\cup\h0{\nu-s+1}\cup\h1{\nu+s'}$
where $1\le s=$ 
$1+c+(n-1)d\le \nu-1$, $0\le c, d \le n-1$,
$s'=d+(n-1)c$.  Note that, if $G$ is well-placed 
with respect to the coarse handle decomposition, then $G$ is 
certainly well-placed with respect to the fine one.

\begin{theorem}\label{theorem 4.2}
Any full subsurface of $T_{n,n}$ is quasipositive.
\end{theorem}

\begin{proof}
Let $S\sub S(\brep{q})$ be a full subsurface.  We may assume 
$S\cap\Bd S(\brep{q})=\varnothing$.  Among all 
graphs $G\sub S(\brep{q})$
with $N_{S(\brep{q})}(G)$ isotopic to $S$ on $S(\brep{q})$, 
such that $G$ is proper,~full,~trivalent,~and 
well-placed \textit{with respect to the coarse handle decomposition}, 
let $G_0$  be such that 
$\skel{G_{\scriptscriptstyle 0}}{1}$ (also %
with respect to the coarse handle decomposition) has the minimal number 
of components.   
Now apply \ref{lemma 3.2}, 
\textit{with respect to the fine handle decomposition}, 
to find a combed graph $G_1$   
isotopic to $G_0$.  
Let $E\sub G_1$  be a comb; then $E$ has no more than 
one endpoint in any   
attaching arc of a $1$-handle of the 
fine handle decomposition, for otherwise \ref{lemma 3.3} 
would contradict the assumed minimality of $G_0$.

We are nearly done.  Of course $N_{S(\brep{q})}(G_1)$ is 
isotopic to $S$ on $S(\brep{q})$.  On the other hand, 
Fig.~\ref{Figure 5} shows how to perform local moves (in the vicinity 
of each coarse $1$-handle) which 
effect an ambient isotopy of $N_{S(\brep{q})}(G_1)$ is 
in $\R^3$ that pushes $N_{S(\brep{q})}(G_1)$ off $S(\brep{q})$ 
and onto a quasipositive braided surface.  
(The fact that the coarse decomposition has a single\linebreak[4]

\pagebreak[4]
\markright{\protect\small Lee Rudolph}
\begin{figure}[t]
\centering
\includegraphics[width=4.5in]{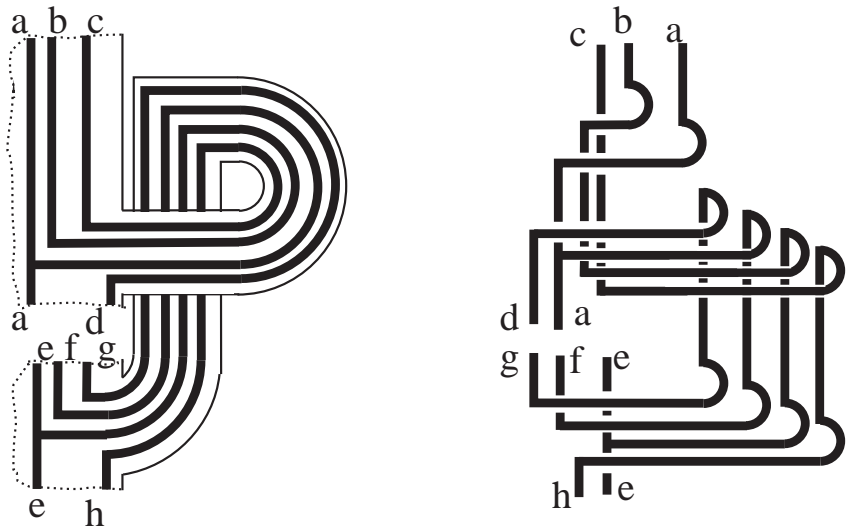}
\caption{The black surface on the left is 
\protect$N_{\scriptscriptstyle{S(\brep{q})}}(G_{\scriptscriptstyle 1})$;
it is isotopic, in $\R^3$, to the black surface on the right, which
\protect\phantom{Fig. 5. The black s}can be braided 
by expanding the thin vertical pieces into fat $0$-handles.\hfill
\label{Figure 5}}
\end{figure}
\noindent
$0$-handle ensures that the local moves do not interfere with each other.  
It was for this reason 
that the coarse decomposition was introduced.) 
\end{proof}

Taken with \S\ref{second section}, \ref{theorem 4.2} completes the 
proof of the Characterization Theorem.

\section{A CONJECTURAL GENERALIZATION}\label{generalization}
Let $d:\C^2\to [0,\infty{[}:(z,w)\mapsto (|z|^2+|w|^2)^{1/2}$.  
A smoothly embedded surface 
$S\sub D^4(r)\isdefinedas d^{\scriptscriptstyle{-1}}([0,r])$, $r>0$, is 
\bydef{ribbon-embedded} if $\Bd S=S\cap \Bd D^4(r)$
and the restriction $d^2|S$ is a Morse function which has no 
local maxima on $\Int S$.  If $S\sub D^4(r)$ 
is ribbon-embedded and 
$r'{\thinspace\scriptstyle\in\thinspace}{]}0,r{[}$ is 
a regular value of $d^2|S$, then $S'\isdefinedas S\cap D^4(r')$ 
is ribbon-embedded in $D^4(r')$, and $S'$ is full on $S$; 
call $(S,S')$ a \bydef{ribbon-embedded pair}.  A surface 
ambient isotopic to a ribbon embedded surface is a \bydef{ribbon}; 
a subsurface $S'\sub S$ of a ribbon is a \bydef{subribbon} 
if the pair $(S,S')$ is ambient isotopic to a ribbon-embedded pair. 

Let  $\Gamma_n(\e)$ be the complex plane curve 
$\{(z,w){\thinspace\scriptstyle\in\thinspace}\C^2\Suchthat
z^n+w^n=\e\}$.  Then, for all sufficiently 
small $\e\ne 0$, the sets $\Gamma_n(\e)\cap D^4(1)$ are 
mutually isotopic ribbon-embedded surfaces.  Let 
${\tilde T}_{n,n}$ denote any one of them.

\begin{conjecture}\label{the conjecture}
A ribbon is quasipositive if and only if, for some $n$, it 
is ambient isotopic to a subribbon of ${\tilde T}_{n,n}$.
\end{conjecture}

Here, by appealing to results in \cite{Rudolph1983b}, 
we can define \bydef{quasipositive ribbon} as follows 
(the original definition involves braids, and ``band representations'' 
more general than the embedded band representations used in this paper).  
Let $D$ denote the bidisk $\{(z,w)\Suchthat |z|\le 1, |w|\le 1\}$,  
$\Bd_1 D\isdefinedas \{(z,w)\Suchthat |z|=1, |w|\le 1\}$.
Let $\Gamma$ be a non-singular complex plane curve which intersects  
$\Bd_1 D$ transversely and the rest of $\Bd D$ not at all.  
Let $\eta: D\to D^4(1)$ be a smoothing (a homeomorphism 
which is a diffeomorphism except along the corner 
torus $\Bd(\Bd_1 D))$.  Then $\eta(\Gamma\cap D)$
is a quasipositive ribbon, and every quasipositive ribbon arises this 
way.

\pagebreak[4]
The following facts are known, cf.\ \cite{Rudolph1983a},
\cite{Rudolph1983b},
\cite{Rudolph1983c}.  
\begin{inparaenum}
\item
One may ``push in'' the interior of a~Seifert surface 
in $S^3=\Bd D^4(1)$ to obtain a ribbon with the same boundary; 
if the Seifert surface is quasipositive, so is the 
ribbon.
\item
In particular, ${\tilde T}_{n,n}$ itself is a quasipositive\break
ribbon, 
for it is the push-in of $T_{n,n}$ (in fact, they are the two 
versions---one in the $4$-disk, the other in the $3$-sphere---of 
the Milnor fiber of $z^n+w^n$).
\item
If $S$ is a quasipositive ribbon, then for some $n$, 
$S$ is a subribbon of ${\tilde T}_{n,n}$.  
(This can be seen as follows.  Realize $S$ by~$\eta(\Gamma\cap D)$~as\break 
above.  Without loss of generality, the completion of 
$\Gamma$ in $\C\P^2\supset\C^2$ is non-singular 
and transverse to the line at infinity.   
Let $f(z,w)$ be the defining polynomial of $\Gamma$.  
For all sufficiently large $R>0$, $\eta(\Gamma\cap D)$ 
is isotopic to a union of components of the ribbon-embedded 
surface $\{(z,w){\thinspace\scriptstyle\in\thinspace}\C^2\Suchthat
f(z,Rw)=0\}\cap D^4(1)$ in $D^4(1)$.  This 
is in turn a subribbon of 
$\{(z,w){\thinspace\scriptstyle\in\thinspace}\C^2\Suchthat f(z,Rw)=0\}
\cap D^4(r)$ in $D^4(r)$ for generic $r\ge 1$.  
Because $\Gamma$ is transverse to the line at infinity, 
$\{(z,w){\thinspace\scriptstyle\in\thinspace}\C^2\Suchthat f(z,Rw)=0\}
\cap D^4(r)$ is isotopic to $\{(z,w){\thinspace\scriptstyle\in\thinspace}\C^2%
\Suchthat z^n+w^n=1\}\cap D^4(r)$ for all sufficiently 
large $r$, where $n$ is the degree of $f$.  
Finally, by rescaling, this last ribbon is isotopic to 
${\tilde T}_{n,n}$.)
\end{inparaenum}

Thus what remains conjectural is that every subribbon 
of ${\tilde T}_{n,n}$ is quasipositive.  If this is true, 
then (up to ambient isotopy) an oriented link $L$ in 
$S^3=\Bd D^4(r)$ bounds a piece of complex plane curve 
in $D^4(r)$ if and only if $L$ has some representation 
as the closure of a quasipositive 
braid.\appnote{Boileau-Orevkov} The 
Characterization Theorem gives some hope that this is so, 
but it seems likely that entirely different techniques 
will be needed for a proof.


\begin{thebibliography}{1}
\setcounter{enumiv}{-1}


\bibitem{Lyon}
\textsc{H.~L\scriptsize{yon}}: 
{Torus knots in the complements of links and surfaces}, 
\textit{Mich. Math. J.} \textbf{27} (1980), 39--46.

\bibitem{Rudolph1983a}
\textsc{L. R\scriptsize{udolph}}: 
\href{http://arxiv.org/abs/math.GT/0411316}%
{Algebraic functions and closed braids}, 
  \textit{Topology} \textbf{22}
  (1983), 191--202. 

\bibitem{Rudolph1983b}
\textsc{L. R\scriptsize{udolph}}:
\href{http://134.76.163.65/servlet/digbib?template=view.html&%
  id=214086&startpage=5&endpage=41&imageset-id=5267}%
{Braided surfaces and {S}eifert ribbons for closed braids},
  \textit{Comment. Math. Helv.} \textbf{58} (1983), 1--37. 

\bibitem{Rudolph1983c}
\textsc{L. R\scriptsize{udolph}}: 
{Constructions of quasipositive knots and links, \textup{I},}
\textit{N{\oe}uds, Tresses, et Singularit\'es} (ed. C. Weber),  
\textit{L'Ens. Math.} (1983), 233--245.

\bibitem{Rudolph1984}
\textsc{L. R\scriptsize{udolph}}: 
{Constructions of quasipositive knots and links, II,} 
\textit{Contemp. Math.}
\textbf{35} (1984), 485--491.

\bibitem{Rudolph1990a}
\textsc{L. R\scriptsize{udolph}}: 
Quasipositive annuli (Constructions of quasipositive knots and links, IV), 
preprint (1990).\appnote{QPIV updated} 


\bibitem{Rudolph1990b}
\textsc{L. R\scriptsize{udolph}}: 
Quasipositive plumbing (Constructions of quasipositive knots and links, V), 
preprint (1990).\appnote{QPV updated} 

\bibitem{Rudolph1989}
\textsc{L. R\scriptsize{udolph}}: 
Quasipositivity and new link invariants, 
\textit{Revista Matem\'atica de la Universidad 
Complutense de Madrid} \textbf{2} (1989), 85--109.

\bibitem{Stallings}
\textsc{J.~R.~Stallings}: {Constructions of fibred knots and links}, 
\textit{Proc. Sympos. Pure Math.} XXXII, Part 2 (Providence: AMS, 1979), 
55--60.   

\end{thebibliography}

\begin{thebibliography}{10}
\setcounter{enumiv}{9}

\bibitem{Boileau-Orevkov}
\textsc{M\scriptsize{ichel} B\scriptsize{oileau}} and 
\textsc{S\scriptsize{tepan} Y\scriptsize{u}. O\scriptsize{revkov}}:
{Quasipositivit\'e d'une courbe analytique dans une boule pseudo-convexe},
\textit{C. R. Acad. Sci. Paris} \textbf{332} (2001), 825--830.
\MR{1836094}

\bibitem{Rudolph1990a updated}
\textsc{L. R\scriptsize{udolph}}: 
\href{http://arxiv.org/abs/math.GT/0112277}%
{Quasipositive annuli (Constructions of quasipositive knots and links, IV)}, 
\textit{J. Knot Theory Ramifications}
\textbf{1} (1992), 451--466.
\MR{1194997}

\bibitem{Rudolph1990b updated}
\textsc{L. R\scriptsize{udolph}}: 
\href{http://www.ams.org/journal-getitem?pii=S0002-9939-98-04407-4}%
{Quasipositive plumbing (Constructions of quasipositive knots and links, V)}, 
\textit{Proc. Amer. Math. Soc.} \textbf{126}, 257--267.
\MR{1452826}

\end{thebibliography}
\renewcommand{\MR}[1]{%
  \href{http://www.ams.org/mathscinet-getitem?mr=#1}{MR#1}
}

{\obeylines\small
\textit{
\noindent Department of Mathematics and Computer Science
\noindent Clark University
\noindent Worcester
\noindent Massachusetts 01610
\noindent U.S.A.}
}

\newpage
\thispagestyle{empty}
\normalsize\normalfont
\section*{ADDENDA}
\hypertarget{addenda}{}
Typographical errors in the original publication 
have been corrected without notice; it is to be hoped that 
no new ones have been introduced. 
The following notes provide updates on various points.

\begin{asparaenum}
%
%
\item\label{Boileau-Orevkov}
Boileau and Orevkov \cite{Boileau-Orevkov} 
have proved that, indeed, up to ambient isotopy
an oriented link $L$ in $S^3$ bounds a piece of complex 
plane curve in $D^4$ if and only if $L$ has some 
representation as the closure of a quasipositive 
braid.  In fact their proof shows that many
subribbons of ${\tilde T}_{n,n}$ are quasipositive;
as far as I am aware, the Conjecture in full generality 
remains open (but is perhaps of somewhat less interest).
\item\label{QPIV updated}
\cite{Rudolph1990a} was published as \cite{Rudolph1990a updated}.
\item\label{QPV updated}
A considerably expanded version of \cite{Rudolph1990b} was 
published as \cite{Rudolph1990b updated}. 

\end{asparaenum}

\renewcommand\refname{Additional References}

\end{document}